\documentclass[leqno,12pt]{article}
\title{Controllability of NLS in the vicinity of solitary wave solutions}
 \author{H.~Lange (Cologne)\\
 H.~Teismann (Acadia)}
\usepackage{amsmath}
\usepackage{amssymb}
\usepackage{amscd}
\usepackage{mathrsfs}
\newcommand{\R}{\mathbb{R}}
\newcommand{\C}{\mathbb{C}}
\newcommand{\N}{\mathbb{N}}

\newcommand{\Z}{\mathbb{Z}}
\newcommand{\cn}{\textrm{cn}}
\newtheorem{thm}{Theorem}

\newtheorem{lem}[thm]{Lemma}
\renewcommand{\Re}{\mathrm{Re}}
\renewcommand{\Im}{\mathrm{Im}}
\newcommand{\muz}{\mu}
\renewcommand{\u}{y}
\renewcommand{\chi}{g}
\newcommand{\g}{u}

\begin{document}

\maketitle

 \begin{abstract}
 Local exact controllability of the 1D NLS  (subject to
 zero boundary conditions) with distributed control is shown  
 to hold in a $H^1$--neighbourhood of the nonlinear ground state.   
 The \textit{Hilbert Uniqueness Method (HUM)}, due to 
 J.-L. Lions, is applied to the linear control problem 
that arises by linearization around the ground state. 
The application of HUM crucially depends on the 
spectral properties of the linearized NLS operator which are given in
detail.   
  \end{abstract} 
\section{Introduction}
The control properties of many PDEs arising in physics
and engineering have been studied extensively. Those investigations 
include exact and/or optimal controllability for the (linear and 
nonlinear) heat, wave, beam and plate equations as well as
 equations of elasticity and Navier Stokes equations, to mention just some of 
the most prominent examples. For Schr\"odinger equations, however,  the 
control theory is markedly 
less developed. For a brief survey on control results 
for (linear and nonlinear) Schr\"odinger equations, see e.g. \cite{ILT06}.
The purpose of this paper is to establish local exact controllability
for the 1D cubic nonlinear Schr\"odinger equation (subject to
zero boundary conditions) with distributed control. 
Specifically, we will consider the following control problem.  
\begin{subequations}
\begin{eqnarray}
 i\u _t &=& -\u _{xx}+f(|\u |^2)\u +\chi _{_{\Omega '}}(x)\g(t,x), 
\quad (x\in (0,1),\, t\in (0,T))          \label{NLSa} \\
 \u (t,0)&=&\u (t,1)=0 \quad (t\in [0,T]) \label{NLSb}\\
 \u (0,x)& =& \u _0(x)\quad (x\in (0,1))  \label{NLSc}\\
 \u (T,x)&=& \u _1(x) \quad (x\in (0,1))  \label{NLSd}
\label{NLS}
\end{eqnarray}
\end{subequations}
where $\chi _{_{\Omega '}}(x)$ denotes the indicator function for 
some fixed, possibly small, open subinterval 
$\Omega '\subset \Omega := (0,1)$ which
represents the spatial region in which the control is applied. 
Given a fixed control ``horizon'' $T>0$ and
initial and target states $\u _0$ and 
$\u _1$, the objective is  to construct a control function 
$\g=\g(t,x)$ that will steer the state $\u $ from $\u _0$ to $\u _1$, i.e. 
the unique solution $\u =\u (t,x)$ of (\ref{NLSa})-(\ref{NLSc}) is to satisfy
(\ref{NLSd}). The control problem (\ref{NLSa})-(\ref{NLSd}) 
was posed  in \cite{ZZ02,ZZ03}. A small-data controllability 
result for periodic boundary conditions is contained in  \cite{ILT02}.\\
In the remainder of this paper, we will concentrate
on the special case of the focusing cubic nonlinear Schr\"odinger 
equation (NLS), i.e.
\[ f(s) = -s .\] 
More general nonlinearities could be treated, but this is not our 
main interest here. So we restrict ourselves to the prototypical cubic
nonlinearity. Our main result states that the control problem 
(\ref{NLSa})-(\ref{NLSd}) is soluble locally within an $H^1$--neighbourhood 
of the (ground-state)  \textit{solitary-wave} (or \textit{soliton}) solution
 $\varphi _\mu (t,x) = e^{i\muz t} \phi _\mu (x)$, if the control time 
$T>0$ is sufficiently large. Here $\phi _\mu (x)$ denotes the  (nonlinear) 
ground state\footnote{Often the time-dependent  solution $\varphi _\mu (t,x)$ 
is referred to as the ground state.}, 
i.e.  the (real and) positive solution 
of the boundary value problem 
\begin{subequations}
\begin{eqnarray}
 -\phi '' &+&\muz  \phi -\phi ^3=0 \quad (x\in \Omega )\label{groundA} 
                                                       \label{groundA1} \\
 \phi (0) &=& \phi (1) =0. \label{groundB} 
                           \label{groundB1} 
 \end{eqnarray}
\end{subequations} 
which is known to exist and to be unique (see Section \ref{Elliptic}). 
Our main results reads. 
\begin{thm} \label{T} Let $\Omega '\subset (0,1)$, $\mu >0$
be given and let $\phi _\mu $ denote the ground state.
 Then there exist $T>0$ and  $\delta >0$ such that, for any 
 $\u _0$, $\u _1\in H^1_0(0,1) := $ \linebreak 
 $\{ v\in H^1(0,1)\mid v(0)=v(1)=0\} $ satisfying 
\[ \| \u _0-\phi _\mu \|_{H^1} <\delta \quad \textrm{and} \quad
 \| \u _1-e^{i\muz T}\phi _\mu \|_{H^1} <\delta  ,\]  
the control problem (\ref{NLSa})-(\ref{NLSd}) has a solution 
 $\g=\g(t,x)$, i.e. there exists a control function 
$\g\in C([0,T];H^1_0(0,1))$ such that
the unique solution $\u \in C([0,T];H^1_0(0,1))$ of 
(\ref{NLSa})-(\ref{NLSc}) satisfies (\ref{NLSd}). 
\end{thm} \label{T1}
The theorem will be proved by applying the implicit function theorem (IFT)
to the nonlinear map 
 $\Psi : [H^1_0(0,1)]^2\times C([0,T];H^1_0(0,1))\to H^1_0(0,1)$, 
defined by  
\begin{equation}
 \Psi (\u _0,\u _1,\g) := \u (T;\u _0,\g) -\u _1, \label{map}
 \end{equation} 
where $t\mapsto \u (t;\u _0,\g)\in  C([0,T];H^1_0(0,1))$ 
denotes the unique solution of (\ref{NLSa})-(\ref{NLSc}). 
To be able to apply the IFT, it will be verified that the linearization 
$\partial _\g\Psi (\phi _\mu ,e^{i\muz T}\phi _\mu,0):
 C([0,T];H^1_0(0,1))\to H^1_0(0,1)$ 
of $\Psi $ at the point $(\phi _\mu ,e^{i\muz T}\phi _\mu,0)$ exists as
a bounded map and possesses a bounded inverse. This amounts to 
showing that the linear PDE that arises by linearizing NLS around the 
stationary solution $\varphi _\mu $ is exactly controllable; this, in turn,
 is done by employing the \textit{Hilbert Uniqueness Method} (HUM) due to
J.-L. Lions \cite{Lio88}. The main difficulty in the application of 
HUM stems from the lack of selfadjointness of the 
linearized operator, which makes determining its 
spectral properties more intricate. 
The analysis reveals that most spectral 
properties known to hold for the 
linearized NLS operator in the \textit{whole-space} case 
\cite{RSS05,CGNT06},
carry over to the \textit{zero-boundary} case considered in this paper
(see Section \ref{Spec})\footnote{Obviously, this statement is not
meant to be applied to those properties which crucially depend on the 
fact that the spatial domain $\Omega =(0,1)$ is bounded, such as 
the absence of a continuous part of the spectrum.}. 
This fact is of independent interest, but it is, to the best of our knowledge,
not available in the literature.  \\[0.1ex]
This paper is organized as follows. In Section 2 we formulate the 
linear control problem that arises by linearization around the 
ground state (system (\ref{Za})-(\ref{Zd})) and state its solvability
(Th. \ref{P}). We also show how Theorem \ref{T1} is derived from this
controllability result (Section \ref{Proof}). 
Section 3 contains the proof of Th. \ref{P}. The proof is based on 
HUM and, thus, hinges on the ``observability estimates'' that are proved in 
Sections \ref{L2ob} ($L^2$-estimate) and 
\ref{H1ob} ($H^1$-estimate). The all-important spectral properties
needed in these proofs are listed in Section \ref{Spec} and proved in 
the Appendix (\ref{Ver}). The Appendix also contains information 
on explicit solution formulas for  (\ref{groundA}),(\ref{groundB}) 
(Section \ref{Elliptic}), a variational  description of the 
ground state (\ref{Var}),  and asymptotic formulas for the 
eigenvalues and eigenfunctions (\ref{Asymp}). 
A number of open problems are  listed in Section \ref{Con}.      
\section{Linearization and proof of Theorem \ref{T1}}
Let $T>0$ be such that property \#\ref{10} in \ref{Spec} is satisfied.
(Such a $T$ exists according to \#\ref{10} in \ref{Ver}.) 
The parameters $T>0$, $\mu >0$ and the interval $\Omega '\subset (0,1)$ 
will be kept fixed in all what follows. We also assume 
 that $\chi := \chi _{_{\Omega'}} $ is a smooth function
satisfying  $\textrm{supp} (\chi )\subset \Omega '$ and 
$0\le \chi  (x)\le 1$.\footnote{This assumption  can be made without
loss of generality. To see this, assume 
that Theorem \ref{T1} is proved for this case. 
If $\chi _{_{\Omega '}}$ is the ``actual'' (non-smooth) 
indicator function  for the interval $\Omega '$,  choose 
an open subinterval $\tilde{\Omega } \subset \Omega '$ 
and a smooth function $\tilde{\chi } $ with 
$\textrm{supp} (\tilde{\chi })\subset \tilde{\Omega }$
 and $0\le \tilde{\chi }(x)\le 1$. Then, by Theorem 
\ref{T1}, there will be a control $\tilde{\g} $ that solves the 
control problem (\ref{NLSa})-(\ref{NLSd}) with $\Omega '$ replaced by 
$\tilde{\Omega }$.
Now $\g (t,x) := \tilde{\chi } (x)\tilde{\g} (t,x)$ will 
be a suitable control for the original problem.} 
\subsection{Linearization of $\Psi $}
We write $\Psi $ as 
\[ \Psi (\u _0,\u _1,\g) = \Phi (\u _0,\g) -\u _1,\] 
where $\Phi :H^1_0(0,1) \times  C([0,T];H^1_0(0,1))\to H^1_0(0,1)$
is the map $\Phi (\u _0,\g) := \u (T;\u _0,\g)$. 
To see that the map $\Phi $ (and hence $\Psi $) 
is well-defined, we need to know that   
the initial value problem (\ref{NLSa})-(\ref{NLSc})
has a unique solution $\u \in C([0,T];H^1_0(0,1))$ for any choice of data
$\u _0\in  H^1_0(0,1)$ and $\g\in C([0,T];H^1_0(0,1))$. 
This is known for the homogeneous NLS 
(i.e. $\g\equiv 0$ in (\ref{NLSa})) in 1D; cf. \cite[Corollary 3.5.2.]{Caz03}.
It is fairly easy to convince oneself that NLS with  an additional 
inhomogeneity  $\tilde{\g} \in C([0,T];H^1_0(0,1))$ (which is given by 
$\tilde{\g} (t,x)= \chi   (x) \g (t,x) $ in our case)
can be treated with the same methods as the the homogeneous equation.
It is also not difficult to verify that the map $\Phi $ (and, by extension,
$\Psi $) is continuous and Fr\'echet differentiable. We omit the 
technical details.  Note that 
\begin{subequations}
\begin{eqnarray}
\Phi (\phi _\mu,0) &=& \varphi _\mu (T) = e^{i\muz T} \phi _\mu  
\quad \Rightarrow \quad \Psi (\phi _\mu ,e^{i\muz T}\phi _\mu ,0) =0
\\
\partial _\g\Psi (\phi _\mu ,e^{i\muz T} \phi _\mu ,0)\cdot h &=&
\partial _\g\Phi (\phi _\mu  ,0)\cdot h \quad 
(\forall h\in  C([0,T];H^1_0(0,1))).
\end{eqnarray}
\end{subequations}
Moreover, the derivative
$\partial _\g\Phi (\phi _\mu ,0):C([0,T];H^1_0(0,1))\to H^1_0(0,1)$
of $\Phi $ at the point $(\phi _\mu ,0)$ 
is given by 
\[ \partial _\g\Phi (\phi _\mu ,0) \cdot h = z(T;h) ,\]
where $t\mapsto z(t)=z(t;h)$ is the solution of IBVP
\begin{subequations}
\begin{eqnarray}
 iz_t &=&-z_{xx} +f(|\varphi _\mu|^2) z 
                +2f'(|\varphi _\mu|^2)\Re (\varphi _\mu\bar{z})\varphi 
 +\chi   (x) h(t,x) \label{Va}\\
 z(t,0)&=&z(t,1)=0  \label{Vb} \\
 z(0,x) &\equiv & 0 \label{Vc}
\end{eqnarray}
\end{subequations}
As usual, the time dependence of the term involving 
$\Re (...)$ is eliminated by the transformation  
\[ \tilde{z} (t) := e^{-i\muz  t} z(t) , 
\quad \tilde{h} (t,x) := e^{-i\muz t} h(t,x) ,\]   
which, for $f(s) =-s$ (cubic focusing nonlinearity), gives the IBVP 
\begin{subequations}
\begin{eqnarray}
 i\tilde{z}_t &=&-\tilde{z}_{xx} +\muz \tilde{z}-\phi _\muz ^2\tilde{z}
 -2 \phi _\muz^2  \Re (\tilde{z})+\chi   (x) \tilde{h}(t,x) \label{VVa}\\
 \tilde{z}(t,0)&=&\tilde{z}(t,1)=0 \label{VVb} \\
 \tilde{z}(0,x) &\equiv & 0 \label{VVc}
\end{eqnarray}
\end{subequations}
We will work with the real $(2\times 2)$-system arising from 
(\ref{VVa})-(\ref{VVc}) by decomposition in real and imaginary parts.
(We will drop the $\mu $ subscripts whenever there is no danger of ambiguity.)
Consider the matrix operator 
\begin{equation}
 L:= \left( \begin{array}{cc} 0 & -\Delta +\muz  -\phi ^2(x) \\
 \Delta -\muz  +3\phi ^2(x) & 0 \end{array} \right) 
 =:  \left( \begin{array}{cc} 0 & L_- \\
 -L_+& 0 \end{array} \right) ;
 \label{L}
 \end{equation}  
this operator is to be understood as an operator in the  Hilbert space 
$H^{-1} \times H^{-1} = [H^1_0(0,1)]^*\times [H^1_0(0,1)]^*$ 
with domain $H^1_0(0,1)\times H^1_0(0,1)=[H^1_0(0,1)]^2$. 
Then eq. (\ref{VVa}) takes the form
 \begin{equation}
  Z_t = LZ + \chi   (x) \left( \begin{array}{r} \Im (\tilde{h}(t,x)) \\ -\Re (\tilde{h}(t,x)) 
 \end{array} \right) , \label{Z} 
 \end{equation} 
 where 
 \[ Z(t,x) = {\Re(\tilde{z}(t,x)) \choose \Im (\tilde{z}(t,x))} .\] 
The corresponding linearized control problem to be considered therefore reads.
\begin{subequations}
\begin{eqnarray}
 Z_t &=& LZ + \chi   (x)H(t,x) \label{Za} \\ 
 Z(t,0) &=& Z(t,1)=\textstyle {0\choose 0} \label{Zb} \\
 Z(0,x) &=& \textstyle {0\choose 0} \label{Zc}\\
 Z(T,x) &=& Z_1(x) \label{Zd}  
 \end{eqnarray}
\end{subequations} 
Our main objective is to prove the following controllability result for 
this linear system. 
\begin{thm} \label{P} The control problem (\ref{Za})-(\ref{Zd}) 
 has a solution $H=H(t,x)$, i.e. for any $Z_1\in [H^1_0(0,1)]^2$ there exists 
 a control function $H\in C([0,T];[H^1_0(0,1)]^2)$ such that the solution 
 $Z\in C([0,T];[H^1_0(0,1)]^2)$ of (\ref{Za})-(\ref{Zc}) satisfies
 (\ref{Zd}). Moreover, the control $H$ satisfies the estimate 
 \begin{equation}
 \sup _{0\le t\le T} \|H(t) \|_{H^1} \le C\|Z_1\|_{H^1} .\label{est}
\end{equation} 
\end{thm}
 Before proceeding with the proof of this theorem, we describe how 
 Theorem \ref{T1} follows from it.  
 \subsection{Proof of Theorem \ref{T1}} \label{Proof}  
As mentioned in the introduction, our plan is to apply the IFT (see, e.g. 
\cite[Theorem 4.B]{Z86}) to the 
map $\Psi :[H^1_0(0,1)]^2\times C([0,T];H^1_0(0,1))\to
H^1_0(0,1)$, $\Psi (\u _0,\u _1,\g):=\Phi (\u _0,\g) -\u _1$. Since
\begin{equation}
 e^{-i\muz  T}\partial _\g\Phi (\phi,0)\cdot e^{i\muz  T}\tilde{h} = z(T;\phi,\tilde{h} ) .
 \label{note}
\end{equation}
the map $\partial _\g \Phi (\phi,0)$ is onto iff
the equation (\ref{VVa})-(\ref{VVc}) is exactly controllable, 
which, in turn, is the case precisely if the  equation 
(\ref{Za})-(\ref{Zc}) is exactly controllable. Thus, by Th. ref{P},
 $\partial _\g \Phi (\phi,0)$ (and hence  
 $\partial _\g \Psi (\phi,e^{i\mu ^2T},0)$) is onto. Moreover, 
\[ \| \partial _\g\Phi (\phi,0)\cdot h\| = \| z(T) \| _{H^1} 
 = \| Z(T) \|_{H^1} \ge C \|H \| _{H^1} =C\| h\|_{H^1},\]
 which shows that the map $\partial _\g \Phi (\phi,0)$ is one-to-one as well 
and that its inverse is bounded. Thus, by the IFT, there exist neighbourhoods 
$Y_0\times Y_1\subset [H^1_0(0,1)]^2$ and $U\subset C([0,T];H^1(0,1))$ of 
$(\phi ,e^{i\muz T} \phi )$  and $0$, respectively, such that
  \[ \forall (\u _0,\u _1) \in Y_0\times Y_1\, \exists ! \g = \g(\u _0,\u _1)\in U : \quad 
  \Psi (\u _0,\u _1,\g) =0 ,\]
  which concludes the proof of the theorem (assuming that Th. \ref{P} has been
 proven), since $\u _1 \in H^1_0(0,1)$ is a reachable state 
iff  there exists a function $\g\in C([0,T];H^1(0,1))$ such that  
$\Psi (\u _0,\u _1,\g)=\Phi (\u _0,\g) -\u _1 =0 $.     
 \hfill $\Box $  
\section{Exact controllability of linearized NLS}
The purpose of this section is to prove Theorem \ref{P}. The proof
will employ the \textit{Hilbert Uniqueness Method (HUM)} 
due to J.-L. Lions \cite{Lio88}.  
We will be using the decomposition $L=L^{(0)}+L^{(1)}$ of $L$, given by
\[ L^{(0)} := \left( \begin{array}{cc} 0 & -\Delta +\muz    \\
 \Delta -\muz    & 0 \end{array} \right) \quad
 \textrm{and} \quad 
 L^{(1)} := \left( \begin{array}{cc} 0 & -\phi ^2(x)  \\
 3\phi ^2(x)  & 0 \end{array} \right) .\] 
We will also need the adjoint operator
\begin{equation}
 L^* = \left( \begin{array}{cc} 0 & \Delta -\muz   +3\phi ^2(x) \\
 -\Delta +\muz -\phi ^2(x) & 0 \end{array} \right) =
  \left( \begin{array}{cc} 0 & -L_+ \\
 L_- & 0 \end{array} \right).
 \label{HUM2}
 \end{equation}  
The set-up is similar to the one in \cite{ILT02}. We use capital letters to 
remind ourselves that we are dealing with two-component functions. 
All equations below are to be understood subject to zero boundary conditions
(i.e. are to be interpreted in the Hilbert space $[H^1_0(0,1)]^2$).
The equation 
\begin{equation}
 \left\{ \begin{array}{ll}
 W_t = LW+g \cdot V \\
 W(T,x)=Z_1(x)  \end{array} \right.
\label{LSb} ,
\end{equation}
is decomposed into two ``semi-homogeneous'' equations:   
\begin{subequations}
\begin{eqnarray}
 W^{(1)}_t &=& LW^{(1)} , \quad 
 W^{(1)}(T,x)=Z_1(x) \label{LSb1} \\ 
 W^{(2)}_t &=& LW^{(2)}+g \cdot V , \quad 
 W^{(2)}(T,x)=0, \quad \textrm{where}   \label{LSb2}  \\
 W&=& W^{(1)} +W^{(2)} \label{W} 
\end{eqnarray} 
\end{subequations} 
The ``HUM operator'' is defined as
\begin{equation}
 SV _0 := -W^{(2)}(0) \label{S} 
\end{equation} 
where $W^{(2)}$ is the solution of (\ref{LSb2}) with $V$ given by 
\begin{equation}
V_t = -L^* V, \quad V(0) = V_0. \label{HUM1} 
\end{equation}  
If $S$ possesses a (bounded) inverse (in the space $[H^1_0(0,1)]^2$), 
then the equation 
\[ SV_0 = W^{(1)}(0) \]
has a unique solution $V_0$, and $W=W^{(1)}+W^{(2)}$ will satisfy (\ref{LSb}) 
with $W(0)=0$. So setting 
\[ Z=W, \quad H =V \]
 will solve the control problem  (\ref{Za})-(\ref{Zd}).
Moreover, the estimate  (\ref{est}) will follow from  the apriori estimates
(\ref{VL2}) and (\ref{VH1}) listed in Section \ref{AP}. 
The crux of the HUM, therefore, consists in showing that 
 $S$ has a bounded inverse. This is done by showing that 
$S$ satisfies an ``observability estimate'' of the form
 \[ \langle   SV_0,V_0\rangle     \,\, \ge C \| V_0\| ^2  \qquad  (\forall V_0\in [H^1_0(0,1)]^2 )\]   
 (\textrm{w.r.t. the appropriate inner product and norm}), which, 
by the Lax-Milgram \linebreak Theorem, implies that $S$ is an isomorphism 
of $[H^1_0(0,1)]^2$.  
\subsection{The $L^2$--observability estimate} \label{L2ob}
 The space $[L^2(0,1)]^2 = L((0,1);\C ^2)$ is equipped with the standard 
inner product  $\langle   .,.\rangle    $, given by
\begin{eqnarray*} \langle   V^{(1)},V^{(2)}\rangle    
 &=&  \left\langle    \left(  \textstyle\begin{array}{c} u^{(1)} \\ v^{(1)} 
 \end{array} \right) ,\left( \begin{array}{c} u^{(2)} \\ v^{(2)} 
 \end{array} \right) \right\rangle      \\
 &=&  \langle    u^{(1)},u^{(2)} \rangle     _{L^2((0,1) ,\C)} + \langle   v^{(1)},v^{(2)}\rangle     _{L^2((0,1) ,\C)}\\  
 &=&  \int _0^{1}u^{(1)}(x)\overline{u^{(2)}(x)} dx
 + \int _0^{1}v^{(1)}(x)\overline{v^{(2)}(x)} dx ;
\end{eqnarray*}
the corresponding norm is $ \|V \|  = \,\sqrt{\langle   V,V\rangle    } $. 
The objective is to show that there exists a constant $C_{HUM}$ such that 
\begin{equation} 
\langle   SV_0,V_0\rangle     \,\, \ge C_{HUM}  \| V_0\| ^2  \qquad  (\forall V_0\in [L^2(0,1)]^2)
\label{HUML2} 
\end{equation}  
 The proof of (\ref{HUML2}) will make use of the following properties
of the operator $L^*$.   
\subsubsection{Spectral decomposition of $-L^*$ and solutions to (\ref{HUM1})}
\label{Spec}
The justification of the properties listed in this section will be deferred  
to the appendix (see \ref{Ver}). 
\begin{enumerate}
 \item \label{1} the spectrum of $-L^*$ consists of eigenvalues only
 \item \label{2} all non-zero eigenvalues are purely imaginary
 \item \label{3} all but finitely many non-zero eigenvalues are simple;
 there are no generalized eigenvectors associated with 
non-zero eigenvalues and each eigenspace is at most two-dimensional.  \\
\textit{Remark.} In the sequel, we will for convenience 
assume that \textit{all} non-zero
 eigenvalues are simple\footnote{There is persuasive numerical evidence 
to support this assumption.}. All subsequent arguments
(in particular, \#\ref{6} below) can easily be 
adapted to accommodate additional (linearly independent) 
eigenvectors which potentially occur for a finite number of eigenvalues.   
 \item \label{4} the multiplicity of the eigenvalue zero is 2; 
 an  eigenvector and a generalized eigenvector, $V_1$ and $W_1$, satisfying 
  \begin{equation}
  -L^*V_1 = \textstyle {0\choose 0} \quad \textrm{and} \quad -L^*W_1=V_1 , \label{jordan}
  \end{equation}
are given by
\begin{equation}
V_1 = {\phi \choose 0}  \quad \textrm{and} \quad 
W_1 ={0 \choose \partial _{\muz   } \phi  } \label{eigenvecs}  
\end{equation} 
(in particular,  $V_1,W_1\in \R^2$). Moreover, these vectors form a basis 
of the generalized null space.  
  \item \label{6} Let $V_2,V_3,\ldots ,\bar{V} _2,\bar{V}_3 ,\ldots $ be 
the eigenvectors corresponding to the simple non-zero eigenvalues 
(see Remark in \#\ref{3} above)
$\lambda _2,\lambda _3,\ldots $, $\bar{\lambda }_2,\bar{\lambda }_3,\ldots $.
Then $V_1,W_1$, $V_2,V_3,\ldots ,\bar{V} _2,\bar{V}_3 ,\ldots $ is a Schauder 
basis for $[L^2(0,1)]^2$ as well as a Bessel sequence, i.e.  
for any $V\in [L^2(0,1)]^2$, there is a unique representation 
\begin{equation} V = c_1V_1+d_1W_1 +\sum _{n\ge 2} [c_nV_n
  +\hat{c}_n\bar{V}_n] \qquad ((c_n,\hat{c} _n) \in \ell ^2)
 \label{Vsum}
 \end{equation}
and there exists a constant $B>0$ (independent of $V$) such that
\begin{equation}
  |c_1|^2+|d_1|^2+\sum _{n\ge 2} [|c_n|^2+|\hat{c}_n|^2]\ge B \|V \| ^2    \label{Bessel}
\end{equation} 
\item \label{7} The (generalized) 
eigenfunctions $V_1,W_1,V_2,V_3,\ldots $ satisfy the uniform bound
\begin{equation}
 m_V:=\min \left\{ \inf _{n\ge 1} \int _\Omega g (x) |V_n(x)|^2 dx 
 , \int _\Omega g (x) |W_1(x)|^2 dx\right\} >0  . \label{inf} 
 \end{equation}
\item  \label{9} Let $V_0\in [L^2(0,1)]^2$ be real, 
$V_0 =  c_1V_1+d_1W_1 +\sum _{n\ge 2} [c_nV_n +\bar{c}_n\bar{V}_n]
= c_1V_1+d_1W_1 +2\Re \sum _{n\ge 2} c_nV_n $ and 
$\lambda _n = i\beta _n$ ($\beta _n \in \R $) .
Then 
\begin{equation} \begin{array}{rcl} 
  V(t) &=& c_1V_1+td_1W_1 +\displaystyle \sum _{n\ge 2} [c_ne^{i\beta _nt}V_n
  +\bar{c}_ne^{-i\beta_nt}\bar{V}_n] \\
 &= & c_1V_1+td_1W_1 +2\Re \displaystyle \sum _{n\ge 2}   c_ne^{i\beta _nt}V_n \label{Vsol}
\end{array}
  \end{equation}
is the unique solution of (\ref{HUM1}). 
 \item \label{10} The sequence 
 \[ 1,\quad t, \quad e^{i\beta _nt}, \quad e^{-i\beta _nt} \qquad (n\ge 2) \]
 is a Riesz-Fischer sequence in $L^2(0,T;\C)$, i.e., there exists a constant 
 $A>0$ such that, for any $\ell ^2$--sequence 
 $a_1,b_1,(a_n,\hat{a}_n)_{n\ge 2}\subset \C $,
 \begin{equation}
  \left. \begin{array}{l} \displaystyle \int _0^T \Big|a_1+b_1t+\sum _{n\ge 2}[a_n e^{i\beta _nt} +\hat{a} _n e^{-i\beta _nt} ]\Big|^2dt \\
 \qquad \qquad \qquad \displaystyle \ge A \Big(|a_1|^2+|b_1|^2+\sum _{n\ge 2} [|a_n|^2+|\hat{a}_n|^2]\Big) \end{array} \right.
  \label{Riesz}
  \end{equation}
\end{enumerate} 
\subsubsection{Proof of (\ref{HUML2})} 
First we show that 
\begin{equation}
 \langle   SV_0,V_0\rangle     \,\,  = \, \int _0^T \langle   gV,V\rangle     dt \label{HUM3} 
\end{equation}
 \textit{Proof of (\ref{HUM3}).} 
Let 
\[ \tilde{S}(t) := \langle   W^{(2)}(t),V(t) \rangle     .\]
Then 
\begin{eqnarray*}
 \frac{d}{dt} \tilde{S}(t) &=&  \langle   W^{(2)}_t,V\rangle    \,+\, \langle   W^{(2)},V_t\rangle     \\
 &\stackrel{(\ref{LSb2}),(\ref{HUM2})}{=} &
  \langle   LW^{(2)},V\rangle     \, + \, \langle    gV,V\rangle     \,-\, \langle   W^{(2)},L^*V\rangle     \,\,= \,\, \langle   gV,V\rangle     
  \end{eqnarray*} 
  and so
\begin{eqnarray*}
  \langle   SV_0,V_0\rangle     &=  & - \, \langle   W^{(2)}(0),V(0)\rangle     \, = \,
 \, \langle   \underbrace{W^{(2)}(T)}_{=0},V(T)\rangle     \, - \, \langle   W^{(2)}(0),V(0)\rangle     \\
  &=&  \tilde{S} (T) - \tilde{S}(0)
  = \int _0^T \frac{d}{dt} \tilde{S} (t) dt = \int _0^T \langle   gV,V\rangle     dt ,
  \end{eqnarray*} 
  which completes the proof of (\ref{HUM3}).   \hfill $\Box $    \\[1ex] 
 Now (\ref{HUML2}) may be verified using standard arguments. 
 The solution $V(t)$ of (\ref{VL2}) has the representation
 \[ V(t,x) = c_1V_1(x)+td_1W_1(x) +\sum _{n\ge 2} 
  [c_ne^{i\beta _nt}V_n(x) +\bar{c}_ne^{-i\beta_nt}\bar{V}_n(x)] \]
 (cf. (\ref{Vsol})). We are going to apply (\ref{Riesz}) with 
 \[ a_1=c_1V_1(x), \quad b_1=d_1W_1(x), \quad a_n = c_nV_n(x), \quad 
  \hat{a} _n=\bar{c} _n\bar{V} _n(x) ;\]
  the result is 
  \begin{eqnarray*} 
\langle   SV_0,V_0\rangle     &\stackrel{(\ref{HUM3})}{=} & \int _0^T \langle   gV,V\rangle     dt =
   \int _\Omega g(x) \int _0^T |V(t,x)|^2 dt dx \\
 &\stackrel{(\ref{Riesz})}{\ge } &
   A  \int _\Omega g(x) \Big(|c_1|^2|V_1(x)|^2+|d_1|^2|W_1(x)|^2
+2\sum _{n\ge 2}
   |c_n|^2|V_n(x)|^2\Big)dx  \\
  &\stackrel{(\ref{inf})}{\ge } & A\, m_V  
 \Big(|c_1|+|d_1|^2 +2\sum _{n\ge 2}
   |c_n|^2\Big) 
  \stackrel{(\ref{Bessel})}{\ge }  A\,B\, m_V  \|V_0\|^2 ,
 \end{eqnarray*} 
which completes the proof of (\ref{HUML2}). \hfill $\Box $ 
\subsection{The $H^1$-observability estimate} \label{H1ob}
Now that the $L^2$--observability estimate 
(\ref{HUML2}) has been established, it is sufficient to show that 
there exist constants $C_1$ and $C_2$ such that 
\begin{equation}
 \langle   (SV_0)_x,{V_0}_x\rangle     \,\, \ge C_1 \| {V_0}_x\| ^2 -C_2 \|V_0\| ^2 \label{HUMH1} 
\end{equation}  
Indeed this will imply 
\begin{eqnarray*}
  \langle   (SV_0)_x,{V_0}_x\rangle     \,\ge  C_1 \| {V_0}_x\| ^2 -\frac{C_2}{C_{HUM}}   \langle   SV_0,V_0\rangle     \\ 
 \Rightarrow \quad 
   \frac{C_2}{C_{HUM}} \, \langle   SV_0,V_0\rangle     \, + \,  \langle   (SV_0)_x,{V_0}_x\rangle     \,\, \ge C_1 \| {V_0}_x\| ^2 .
 \end{eqnarray*} 
 The left-hand side is equivalent to the ``natural'' $H^1$ inner product
 \[ \langle   .,.\rangle     + \langle   (.)_x,(.)_x\rangle     .\] 
Before presenting the proof of (\ref{HUMH1}) we list the apriori estimates 
required. 
 \subsubsection{Apriori estimates} \label{AP}
 We are going to need apriori estimates for the various functions
 involved. The proofs are standard fare and will be omitted.
 \begin{subequations}  
 \begin{eqnarray}
 \| V(t) \|&\le & C \|V_0\|  \label{VL2} \\
 \| V_x(t)\|  &\le & C(\|V_0\|+\|{V_0}_x\| )\label{VH1} \\
 \|W^{(2)}(t)\|  &\le & C\|V_0\| \label{W2L2} \\ 
  \| W^{(2)}_x(t)\| & \le  & C(\| V_0\| +\| {V_0}_x\| )\label{W2H1}
\end{eqnarray}
\end{subequations}
We will also need the equations for 
$V_x$ and $W^{(2)}_x$, which are given by 
\begin{subequations}
\begin{eqnarray}
 {W^{(2)}_x}_t &=& LW^{(2)}_x +L^{(1)}_xW^{(2)}+g_x V +gV_x , \quad 
W^{(2)}_x(T)=0  
\label{LSb2'} \\
{V_x}_t &=& -L^* V_x -[L^{(1)}_x]^*V, \quad V_x(0) = {V_0}_x. \label{HUM1'} ,
\end{eqnarray}    
\end{subequations}         
\subsubsection{Proof of (\ref{HUMH1})}   
We first prove an ``$H^1$ analogy'' of the identity (\ref{HUM3}):
\begin{equation}
 \left\{ \begin{array}{rcl} \langle   (SV_0)_x,{V_0}_x\rangle     &=&  \displaystyle \int _0^T \Big[ \, \langle   L^{(1)}_x W^{(2)},V_x\rangle     \, - \, \langle   L^{(1)}_x W^{(2)}_x,V\rangle    \, \Big] \, dt \\
                   & &    \displaystyle + \int _0^T \langle   g_xV,V_x\rangle     \, dt + \int _0^T \langle   gV_x,V_x\rangle     dt \end{array} \right.  
\label{HUM4}
\end{equation}
\textit{Proof of (\ref{HUM4}).}  Let
\[ \tilde{\tilde{S}} (t) := \, \langle    W^{(2)}_x(t),V_x(t)\rangle     .\]  
Then 
\begin{eqnarray*} \frac{d}{dt} \tilde{\tilde{S}} (t) &=&
 \langle    {W^{(2)}_x}_t,V_x\rangle     \, + \, \langle    W^{(2)}_x,{V_x}_t\rangle     \\
 &=&  \langle   LW^{(2)}_x,V_x\rangle     \,+\, \langle   L^{(1)}_xW^{(2)},V_x\rangle     \,+\, \langle   g_x V,V_x\rangle     \, +\, \langle   gV_x,V_x\rangle     \\
 & &  -\, \langle   LW^{(2)}_x,V_x\rangle     \, -\, \langle   L^{(1)}_xW^{(2)}_x,V\rangle     \\
 &=& \langle   L^{(1)}_xW^{(2)},V_x\rangle      \, -\, \langle   L^{(1)}_xW^{(2)}_x,V\rangle    \,+\, \langle   g_x V,V_x\rangle     \, +\, \langle   gV_x,V_x\rangle     ,
 \end{eqnarray*}     
 which implies (\ref{HUM4}). \hfill $\Box $ \\[1ex]
 Using the abbreviations
 \begin{subequations}
 \begin{eqnarray}
  I_1& :=&  \int _0^T \Big[ \, \langle   L^{(1)}_x W^{(2)},V_x\rangle     \, - \, \langle   L^{(1)}_x W^{(2)}_x,V\rangle    \, \Big] \, dt  \label{I1}\\
  I_2 &:=&   \int _0^T \langle   g_xV,V_x\rangle     dt \label{I2} \\
  I_3 &:=&  \int _0^T \langle   gV_x,V_x\rangle     dt \label{I3}
 \end{eqnarray}
 \end{subequations}
 we estimate 
 \begin{equation}
   \langle   (SV_0)_x,{V_0}_x\rangle     \,\,\stackrel{(\ref{HUM4})}{\ge } \,I_3 -|I_1|-|I_2| .\label{HUM66}
   \end{equation}  
 Now,
 \begin{eqnarray}
  |I_1| &\le & 3\|(\phi ^2)_x\|_\infty \int _0^T(\| W^{(2)} \| \cdot \|V_x\| + \| W^{(2)}_x\| \cdot \|V\| ) dt \nonumber\\
   &\stackrel{(\ref{VL2})-(\ref{W2H1})}{\le } &  C\|V_0\| \cdot (\|V_0\| +\|{V_0}_x\|) \nonumber \\
   &\stackrel{Young}{\le } & \varepsilon \|{V_0}_x\|^2 +C_\varepsilon \|V_0\|^2 \label{HUM8a}.
 \end{eqnarray}   
 The estimate for $I_2$ is even simpler.
 \begin{eqnarray}
  |I_2| &\le & \|g_x\|_\infty \int _0^T\| V \| \cdot \|V_x\|dt   
   \stackrel{(\ref{VL2}),(\ref{VH1})}{\le }   C\|V_0\| \cdot (\|V_0\| +\|{V_0}_x\|)\nonumber \\
   &\stackrel{Young}{\le } & \varepsilon \|{V_0}_x\|^2 +C_\varepsilon \|V_0\|^2 \label{HUM8b} 
 \end{eqnarray}  
 Finally, we need to estimate $I_3$ from below. To do this, we write the 
solution of (\ref{HUM1'}) in Duhamel form
 \begin{equation}
  V_x(t) = e^{-L^*t}{V_0}_x - \int _0^t e^{-L^*(t-s)} [L^{(1)}_x]^*V(s) ds, \label{Vx}
 \end{equation}
 where $e^{-L^*t}$ denotes the semi-group corresponding to the equation 
(\ref{HUM1}).  If we set
 \[ V^{(0)} (t) := e^{-L^*t}{V_0}_x ,\]
 we get from (\ref{HUML2}) and  (\ref{HUM3}) that  
 \begin{equation}
  \int _0^T \int _\Omega g(x) |V^{(0)}(t,x)|^2 dxdt=  \int _0^T\langle   gV^{(0)},V^{(0)}\rangle     dt \ge C_{HUM} \| {V_0}_x\| ^2  
  \label{HUM6}
 \end{equation}  
 (here we mean by $|.|$ the norm of $\C ^2$, i.e. 
$\left| u\choose v \right|^2 = |u|^2+|v|^2 $.)  Now,
 \begin{eqnarray*}
  |I_3| &\stackrel{(\ref{Vx})}{\ge }& \int _0^T \int _\Omega g(x) |V^{(0)}(t,x)|^2 dx dt- \cdots \\
  & &  \quad \quad \qquad \quad \quad - \int _0^T \int _\Omega g(x) \left|\int _0^t e^{-L^*(t-s)} [L^{(1)}_x]^*V(s,x) ds \right|^2 dx dt  \\
  &\stackrel{(\ref{HUM6})}{\ge }&   C_{HUM} \| {V_0}_x\|^2 - \int _0^T \int _\Omega g(x) \left|\int _0^t e^{-L^*(t-s)} [L^{(1)}_x]^*V(s,x) ds \right|^2 dx dt 
 \end{eqnarray*}
 To finish the proof, we need to estimate the integral term. 
 \begin{eqnarray*}
  \int _0^T \int _\Omega g(x) \left|... \right|^2 dx dt   
  &\le &   \int _0^T \int _0^T\int _\Omega  \left|e^{-L^*(t-s)} [L^{(1)}_x]^*V(s,x)  \right|^2 dx ds dt \\
  &=& \int _0^T \int _0^T\left\|e^{-L^*(t-s)} [L^{(1)}_x]^*V(s,.)  \right\|^2  ds dt \\
  &\stackrel{(\ref{VL2})}{\le } &  C\int _0^T \int _0^T\left\|[L^{(1)}_x]^*V(s,.)  \right\|^2  ds dt  \\
 &\le & C \| (\phi ^2)_x\|^2 \int _0^T  \| V(s) \|^2 ds 
 \stackrel{(\ref{VL2})}{\le } C \|V_0\|^2     
 \end{eqnarray*}
 Thus, 
 \begin{equation}
  I_3 \ge C_{HUM}\|{V_0}_x\|^2 - C\|V_0\|^2 \label{HUM7}
  \end{equation}
 and so
  \[ \langle   (SV_0)_x,{V_0}_x\rangle     \,\,\stackrel{(\ref{HUM66})}{\ge } \,I_3 -|I_1|-|I_2|
   \stackrel{(\ref{HUM8a}),(\ref{HUM8b}),(\ref{HUM7})}{\ge } 
   (C_{HUM}-2\varepsilon ) \|{V_0}_x\|^2 - C_\varepsilon \|V_0\|^2   ,\]
  which will conclude the proof of (\ref{HUMH1}) if $\varepsilon $ is chosen 
sufficiently small. \hfill $\Box $   
\section{Concluding remarks} \label{Con}
There are a number of 
 modifications/generalizations of the control problem 
(\ref{NLSa})-(\ref{NLSd}) which are of interest. Those include \vspace{-0ex}  
\begin{itemize}
\item space dimension $>1$;\vspace{-1ex}   
 \item other boundary conditions such as  periodic boundary 
 conditions;   \vspace{-1ex}  
\item  the ground state $\phi _\mu $ in the assumption of 
Theorem \ref{T1} may be replaced with some excited state
 (see \ref{Elliptic});  \vspace{-1ex}  
 \item zero (or ``box'')
 boundary conditions may be interpreted as an
infinite potential well, which one may want to  
replace with other, typically confining, potentials, 
such as the harmonic-oscillator potential (in this case, one would 
reasonably set $\Omega = (-\infty ,\infty )$.)
\end{itemize}
Any of these modifications 
will obviously change some of the spectral properties in \ref{Spec}.
To see whether a controllability result as in Theorem \ref{T1}
can be proved under these modified assumptions as well, 
the ramifications for the application of HUM will have to be carefully 
examined. 
\mbox{}\\
\begin{center}
{\sc \Large \bf Appendix} 
\end{center}
\appendix
\section{Bound states}
We define \textit{bound states} as real solutions of the BVP 
(\ref{groundA1}),(\ref{groundB1}). 
\subsection{Bound states in terms of elliptic functions}
\label{Elliptic}
It is well-known that explicit formulas for the 
solutions of (\ref{groundA1}),(\ref{groundB1})
are available in terms of Jacobian elliptic functions; see, e.g. \cite{CCR00}. 
If $j\in \{0,1,2,\ldots \} $, then $\phi _j(x)$ will denote 
the (real-valued) solution of (\ref{groundA1}),(\ref{groundB1}) 
which possesses precisely $j$ zeros (``nodes'') within the interval
$(0,1)$.  The node-less solution $\phi :=\phi _0$ is referred to as 
the \textit{ground state}; the solutions $\phi _j$ ($j\ge 1$) with one 
or multiple nodes are called \textit{excited states.}  
To find an explicit solution formula for $\phi _j$, we first solve the 
equation 
\begin{equation}
 \muz  = 4(j+1)^2(2k^2-1)K(k)^2 \label{k}
\end{equation}  
for $k$, where $K(k)$ denotes the complete elliptic integral 
of the first kind (see, e.g. \cite{AS65}). Note that,
since $K(k)$ is a strictly increasing continuous function of $k\in [0,1)$ 
satisfying $\lim _{k\to 1^-} K(k) =\infty $, equation (\ref{k}) has 
exactly one solution $k=k_j(\mu )$ for any choice of parameters $\mu \ge 0$ 
and $j\in \{ 0,1,2,\ldots \} $. Moreover, 
the function $k_j: [0,\infty )\to [\frac{1}{\sqrt{2} } ,1)$ 
is continuous and strictly  increasing as well, and satisfies 
$\lim _{s\to \infty} k_j(s) = 1$.
Now the solution $\phi _j$ of (\ref{groundA1}),(\ref{groundB1}) is given by 
($k=k_j(\mu )$)
\begin{equation}
\begin{array}{rcl}
  \phi _j(x) &=&\displaystyle \frac{\sqrt{2\mu} \, k }{\sqrt{2k^2-1}}\,\, \cn \left(
 \frac{\sqrt{\mu} (x-\frac{1}{2} )}{\sqrt{2k^2-1}} +[j]_2K(k),k\right) \\
 &=& 2\sqrt{2} (j+1) k K(k) \,\, \cn \Big(
 2(j+1) K(k) {\textstyle (x-\frac{1}{2} )} +[j]_2K(k),k\Big) ,
\label{phij}
\end{array}
\end{equation}  
where $[j]_2 := j$ mod $2$. The following properties are readily proved. 
\begin{lem} \label{Lell}
Let $\psi \in H^1_0(0,1)$ be a non-trivial weak solution of (\ref{groundA1}),
(\ref{groundB1}), i.e. 
\[ \int _0^1 [\psi _xv_x +\mu \psi v -\psi ^3v ]\, dx =0 \qquad 
(\forall v\in H^1_0(0,1)) .\] 
Then 
\begin{enumerate}
\item[(i)] $\psi \in C^\infty (0,1)\cap C[0,1] $, i.e. $\psi $ is a classical
solution of (\ref{groundA1}),
(\ref{groundB1}). 
\item[(ii)] There exists  $j\in \{0,1,2,\ldots \} $ such that 
 $\psi $ has exactly $j$ zeros in $(0,1) $.
\item[(iii)] $\psi (x) = \phi _j(x)$, where $\phi _j$ is given by (\ref{phij}).
\end{enumerate}
\end{lem}
We also need the following convexity property.  
\begin{lem} \label{L5}
The function $\mu \mapsto  \| \phi \| _2^2$ 
is an increasing function of $\muz $. 
\end{lem}
\textit{Proof.} This follows from the identity
\begin{equation}
\|\phi \|_2^2 = 4 k^2K(k)
 \int _{-K(k)}^{K(k)} \cn ^2(y,k) dy 
\label{phiL2}
\end{equation}
 and the fact that 
$k=k(\mu )$ is a strictly increasing function of $\mu $. \hfill $\Box $
\subsection{Variational description of the ground state} \label{Var}
There are various variational descriptions of the ground state in the 
whole-space case; see e.g. \cite[4.2]{SS99} and \cite[8.1]{Caz03}; 
in the zero-boundary case we find the formulation presented in Lemma \ref{Lvar}
below to be a useful one. Let $\mu _1 $ be the smallest eigenvalue
of the 1D Laplacian $-\frac{d^2}{dx^2 } $ on $(0,1)$ 
with zero boundary 
conditions, i.e. $\mu _1= \pi ^2$. 
\begin{lem} \label{Lvar}
Let $\mu >-\mu_1$. 
Then there exists a positive minimizer 
\[ w\in 
\{ u\in H^1_0(0,1)\mid \int _0^1 u^4 dx =1 \} =:M \]  
for the constrained minimization problem 
\[ \inf _{u\in M } \,  \frac{1}{2} 
\int _0^1 [{u_x}^2 +\mu u^2 ] dx ,\]
and the unique positive solution  $\phi $
of (\ref{groundA1}),(\ref{groundB1}) is given by 
$ \phi (x) = \lambda ^{1/2} w(x) $, 
for some suitable $\lambda >0 $ (Lagrange multiplier).
\end{lem} 
\textit{Proof.} See e.g.  \cite[Theorem 2.1 and proof]{S00}. 
\hfill $\Box $
\section{Spectral properties of $-L^*$}
\subsection{Properties of $L_+$ and $L_-$}  
The operators
\begin{subequations}
 \begin{eqnarray}
 L_-u &= & \textstyle \left[-\frac{d^2}{dx^2}+\muz  -\phi ^2(x)\right] u, \quad 
 u(0)=u(1) =0 \\
  L_+u &= & \textstyle\left[-\frac{d^2}{dx^2}+\muz  -3\phi ^2(x)\right] u, \quad 
 u(0)=u(1) =0
\end{eqnarray} 
\end{subequations} 
 are regular Sturm-Liouville (SL) 
 operators, so we can make use of the SL theory. 
\begin{lem} \label{L1} 
 \begin{enumerate}
 \item[(i)] $\textrm{ker} (L_-) = \textrm{span}\{\phi \} $
 \item[(ii)] $L_-|_{[span (\phi )]^\perp}> 0$  
\end{enumerate}  
\end{lem}
 \textit{Proof. \underline(i)} Clearly, $L_-\phi =0$, 
so $\phi $ is eigenfunction for
$L_-$ with eigenvalue $\lambda =0 $. Since all eigenspaces 
are one-dimensional, the assertion follows. \\[1ex]
\textit{\underline{(ii)}.}  Since $\phi $ has no zeros in 
$(0,1)$, it 
is an eigenfunction for the smallest eigenvalue; hence 
$\lambda _1 =0$ and $\lambda _n >0$ ($\forall n\ge 2$),
which implies assertion (ii).  \hfill $\Box $  \\[1ex]
\textbf{Remark.}
If $\phi $ is not the ground state but an excited state
with 
$j\in \N \setminus \{0\}$ nodes,
 then the operator $L_-$ will have exactly $j$ negative eigenvalues.
\begin{lem} \label{L2}  The operator $L_+$
has exactly one negative eigenvalue and all the other eigenvalues
are positive. In particular,  
 $\textrm{ker} (L_+) = \{0 \} $. 
\end{lem} 
\textit{Proof.} The proof proceeds in four steps. \\
{\sc Step 1.} $L_+$ possess \textit{at least one} negative 
eigenvalue. This follows from  
\[ \langle    L_+\phi ,\phi \rangle      \, = \, \langle    L_-\phi ,\phi \rangle      \, -\, 2\| \phi \| _4^4 
 = -\, 2\| \phi \| _4^4 <0 \]
 and the minimax principle. \\
{\sc Step 2.} $\langle    L_+\eta ,\eta \rangle      \, \ge 0$ for all $\eta \in
[\textrm{span} (\phi ^3)]^\perp $. This is a slight adaptation of the 
arguments in \cite[Section13]{RSS05}; 
we use similar notation. Let the functionals $J[u]$ and $W[u]$ be defined by 
\[ J[u] := \frac{1}{2} \int _0^1 [{u_x}^2+\mu u^2]dx \quad
\textrm{and} \quad W[u]:= \frac{1}{4} \int _0^1 u^4 dx \]
and let $w\in M=\{u\in H^1_0(0,1)\mid W[u]=1\} $ be a positive 
minimizer of the constrained minimization  problem $\inf _{u\in M} J[u] $. 
By Lemma \ref{Lvar}, the ground state $\phi $ is given by 
$\phi (x) = \lambda ^{1/2} w(x)$ for some positive constant $\lambda $, which 
arises as a Lagrange multiplier. Now let 
$\eta \in  [\textrm{span} (\phi ^3)]^\perp $; we write $\eta $ in the form 
$\eta =\dot{w} :=\frac{\partial }{\partial z} |_{z=0} w(.,z)$ 
where $z\mapsto w(.,z)$ is a smooth curve in $H^1_0(0,1)$ such that 
$W[w(.,z)]=1$ (i.e. $w(.,z)\in M$) for all $z$, and $w(.,0)=w$.   
Since $w$ is a minimizer, we have 
\begin{equation}
 0= \frac{d}{dz} \Big|_{z=0} J[w(.,z)]  
= \int _0^1 [w_x\dot{w}_x  +\mu w\dot{w}] dx \label{RSS1}
 \end{equation}
and 
\begin{eqnarray} 0&\le & \frac{d^2}{dz^2} 
\Big|_{z=0} J[w(.,z)]  =
 \int _0^1 [{\dot{w}_x}^2 + w_x\ddot{w}_x  +\mu ({\dot{w}}^2 +w\ddot{w} )] dx
 \nonumber \\
 &=& 
 -\int _0^1 [(\dot{w}_{xx} -\mu {\dot{w}}){\dot{w}} + (w_{xx}  -\mu w)\ddot{w} ] dx
   \label{RSS2}
\end{eqnarray} 
where $\ddot{w} := \frac{\partial ^2}{\partial z^2} |_{z=0} w(.,z)$ Moreover, from the constraint $W[w(.,z)]\equiv 1$
 we get
\begin{equation}
 0 = \frac{d}{dz} \Big|_{z=0} W[w(.,z)]
 = \int _0^1 w^3\dot{w} dx \label{RSS3}
\end{equation}
and 
\begin{equation}
 0 = \frac{d^2}{dz^2} \Big|_{z=0} W[w(.,z)]
 = \int _0^1 [3w^2{\dot{w}}^2+w^3\ddot{w}]  dx \label{RSS4}.
\end{equation}
Also, the Lagrange-multiplier rule implies
\begin{eqnarray} 
\int _0^1 [w_{x}v_x +\mu wv]\, dx &=& \lambda  \int _0^1 w^3 v \, dx 
\quad (\forall v\in H^1_0(0,1))  \label{RSS8} \\
 \stackrel{v=\ddot{w} }{\Rightarrow } \quad 
\int _0^1 [-w_{xx} +\mu w]\ddot{w} \, dx &=& \lambda  \int _0^1 w^3 \ddot{w} dx
\stackrel{(\ref{RSS4})}{=} -\lambda \int _0^1 3w^2{\dot{w}}^2  dx 
\label{RSS5} 
\end{eqnarray}
Inserting this into (\ref{RSS2}) gives 
\begin{eqnarray*} 0&\le & 
 \int _0^1 [-\dot{w}_{xx} +\mu {\dot{w}} -\lambda 3w^2\dot{w} ]\dot{w} dx \\
 &\stackrel{\phi = \lambda ^{1/2}w,\, \eta =\dot{w}}{=} &
 \int _0^1 [-\eta_{xx} +\mu \eta -3\phi ^2\eta ]\eta  \, dx
 = \, \langle    L_+\eta ,\eta \rangle       
\end{eqnarray*} 
(Note that choosing $v=w$ in (\ref{RSS8}) also yields 
$\int _0^1 [w_{x}^2 +\mu w^2]\, dx = \lambda  \int _0^1 w^4 \, dx $
implying that $\lambda >0$, since $\mu >-\mu _1$.) \\
{\sc Step 3.} The second eigenvalue, $\lambda _2$, is non-negative. 
This can be shown by repeating word-by-word the proof in 
\cite[page 58]{RSS05} if $\Psi $ is defined as $\phi ^3$ and 
$\Pi $ is interpreted as the orthogonal projection onto
the subspace $[\textrm{span} (\phi ^3)]^\perp $.  \\
{\sc Step 4.} $\lambda _2>0$. Assume that $\lambda _2=0$ and let 
$v$ denote an eigenfunction for $\lambda _2$. By the symmetry of 
$\phi (x)$, we may assume w.l.o.g. that $v$ is either odd (i.e. $v(1-x) =
-v(x)$) or even (i.e. $v(1-x) = v(x)$). In the first case, we have 
$v(\frac{1}{2} )=0$ and $v$ coincides with a constant multiple of
$\phi '$ by ODE uniqueness. However, this is impossible, since 
$\phi '$ does not satisfy zero boundary conditions. We are therefore 
left with the second case ($v$ even). Since $v$ is an eigenfunction 
for the second eigenvalue, it has precisely one zero in $(0,1)$ by 
standard SL theory. By symmetry this zero must occur at $x=\frac{1}{2} $, 
which is impossible as we saw above.       \hfill $\Box $
\subsection{Eigenvalues and eigenfunctions for $n\to \infty $}
\label{Asymp}
The properties \ref{6} and \ref{7} in Section \ref{Spec} are based on 
asymptotic ($n\to \infty $) formulas for the eigenvalues 
$\lambda _n,\bar{\lambda }_n$ and eigenfunctions 
$V_n(x)$, $\bar{V} _n(x)$ of $-L^*$, which are given in 
Lemmas \ref{LK} and \ref{Lkam} below. 
We are going to make use of the fact (proved in \ref{Ver} \#\ref{2} below) 
that all non-zero eigenvalues are purely imaginary, i.e. we write 
\begin{equation}
 \lambda _n = i\beta _n, \quad \bar{\lambda } _n = -i\beta _n , \qquad
\beta _n\in \R , \, \beta _n>0.\label{beta}
\end{equation}   
Moreover, it will be convenient to employ a similarity transformation
\cite[(12.15)]{RSS05}: Let
\[ J:= \left( \begin{array}{cr} 1 & 1 \\
 i & -i \end{array} \right)  .\]
Then 
\begin{equation} -iL^* = J\cdot  \left[ \left( \begin{array}{cc} \Delta -\muz   &  0 \\
 0 & -\Delta +\muz    \end{array} \right) + \phi ^2 \left( \begin{array}{rr} 2& -1\\
 1 & -2  \end{array} \right)  \right] \cdot J^{-1} .\label{J}
 \end{equation} 
and so
\[ \textrm{spec}(-L^*) = i\textrm{spec}(M) ,\quad 
 V_n = JW_n^+, \quad \bar{V} _n=JW_n^- \] 
where $(\pm \beta _n,W_n^\pm )$ are the eigenpairs for the operator
\[  M:= \left( \begin{array}{cc} \Delta -\muz   &  0 \\
 0 & -\Delta +\muz    \end{array} \right) + \phi ^2 \left( \begin{array}{rr} 2& -1\\
 1 & -2  \end{array} \right) =: M^{(0)} + M^{(1)} ,\]
i.e. 
 \begin{equation} M W_n^\pm = \pm \beta _nW_n^\pm . \label{M3}
\end{equation}  
Writing $\beta =\pm \beta _n$, $W=W_n^\pm ={u\choose v}$, the characteristic 
equation (\ref{M3}) is equivalent to the BVP
\begin{subequations}
\begin{eqnarray}
 u'' -(\muz  +\beta ) u &=& -\phi ^2 (2u-v), \quad \,u(0)=u(1)=0 \label{a}\\
 v'' -(\muz  -\beta ) v &=& \phi ^2 (u-2v), \qquad v(0)=v(1) =0 \label{b} .
\end{eqnarray}   
\end{subequations}
Note that 
\begin{equation}
 W_n^+ = {u\choose v} \quad \iff \quad W_n^-={v\choose u} .\label{uv}
\end{equation}     
\begin{lem} \label{LK} 
\begin{enumerate}
\item[(i)] The operator $M$ (and hence $-L^*$) is a spectral operator.
More precisely, the collection of eigenvectors and generalized eigenvectors 
for $M$ forms a Schauder basis for $L^2(0,1;\C )$ and all eigenvalues 
with sufficiently large indices $n$ are simple.  
\item[(ii)] There exists an index $n_0\in\N $ and a 
constant $C>0$ such that  
   \[ |(n^2\pi ^2+\muz ) -\beta _n | \le C  \]
for all $n\ge n_0$.   
\end{enumerate} 
\end{lem} 
\textit{Proof.} The operator $M^{(0)}$ is self-adjoint; its 
eigenvalues and eigenfunctions are given by 
\[ {\beta _n^{(0)}}^\pm  = \pm (n^2\pi^2+\muz ), \quad
 {W_n^{(0)}}^+ =\sin (n\pi x) {0\choose 1} , \quad 
 {W_n^{(0)}}^- =\sin (n\pi x) {1\choose 0}  \]
(where  $M ^{(0)}{W_n^{(0)}}^\pm = {\beta _n^{(0)}}^\pm {W_n^{(0)}}^\pm  $). 
Thus, the operator $M=M^{(0)}+M^{(1)}$ ``is'' a bounded perturbation of
a self-adjoint operator whose spectrum consists of 
simple eigenvalues only. The assertion now follows from 
\cite[Th. 4.15.a]{K66}; see also
\cite{C68}. \hfill $\Box $
\begin{lem} \label{Lkam}
 There exist an index $n_0\in \N $ and a constant $C>0$ such that  
$\beta _n>\muz $ and   
\begin{subequations}
\begin{eqnarray}
 \left| W_n^+(x) - \sin \left(\sqrt{\beta _n-\muz }\,\, x\right) {0\choose 1} \right| &\le & 
\frac{C}{n}  \label{Wna}\\
 \left| W_n^-(x) - \sin \left(\sqrt{\beta _n-\muz }\,\, x\right) {1\choose 0} \right| &\le & 
\frac{C}{n} \label{Wnb}
\end{eqnarray}
for all $n\ge n_0$.
\end{subequations}
 \end{lem}    
\textbf{Remark.} The assertion of the lemma may be expressed in a more 
intuitive, if slightly informal, manner by means of the asymptotic formulas   
\begin{eqnarray*}
 W_n^+(x) &=& \sin \left(\omega _n\, x\right) \textstyle {0\choose 1} +\mathcal{O} (n^{-1}) , \quad
 n\to \infty \\ 
W_n^-(x) &=& \sin \left(\omega _n\, x\right) \textstyle {1\choose 0} +\mathcal{O} (n^{-1}) , \quad
 n\to \infty .
\end{eqnarray*} 
(uniformly in $x$ and $n$) where 
 $\omega _n := \sqrt{\beta _n-\muz } \sim \pi n $, as $n\to \infty $,
 by Lemma \ref{LK} (ii).  \\[1ex]
\textit{Proof.} Clearly, by Lemma \ref{LK} (ii), 
$\lim _{n\to \infty } \beta _n =\infty $;
so we may assume that $\beta _n-\muz  >0$. Let 
$\omega _n^\pm :=\sqrt{\beta _n\pm \muz } $.  
Because of (\ref{uv}) it is sufficient to consider $W_n^+$; write
$W_n^+(x) ={u(x)\choose v(x)}$. 
Viewing the R.H.S.'s of (\ref{a}),(\ref{b}) as inhomogeneities, 
we write the system as 
\begin{subequations}
\begin{eqnarray}
 u'' -[\omega _n^+]^2 u &=& f(x), \quad u(0)=u(1)=0 \label{af}\\
 v'' +[\omega _n^-]^2 v &=& g(x), \quad v(0)=v(1) =0 \label{bg} .
\end{eqnarray}   
where 
\begin{equation}
 f(x) = \phi ^2(x) [v(x)-2u(x)] \quad \textrm{and} \quad 
 g(x) = \phi ^2(x) [u(x)-2v(x)].\label{fg}
\end{equation}
\end{subequations}
Note that the homogeneous BVP  associated with eq. (\ref{bg}) has the 
solution $\sin(\omega _n^-x)$, while the  homogeneous BVP 
 associated with eq. (\ref{af}) does not have any non-trivial solution, 
which implies that there is a 
Green's function, 
$\Gamma _{\omega _n^+}(x,\xi )$, associated with eq. (\ref{af}).     
Utilizing this Green's function and the Duhamel Principle,
solutions to (\ref{af}),(\ref{bg}) may be written as 
\begin{subequations}
\begin{eqnarray}
 u(x) &=&  \int _0^1 \Gamma _{\omega _n^+}
  (x,\xi) f(\xi )d\xi \label{C1}\\  
 v(x) &=&  c\sin (\omega _n^-x) +\frac{1}{\omega _n^-} \int _0^x
  \sin (\omega _n^-(x-\xi )) g(\xi ) d\xi \quad (c\in \R )\label{D1}  
\end{eqnarray}
\end{subequations} 
By the linearity of the system (\ref{a}),(\ref{b}), we may assume that
 $c=1$. Thus, 
 \begin{subequations}
\begin{eqnarray}
 u(x) &=&  \int _0^1 \Gamma _{\omega _n^+}(x,\xi) \phi ^2(\xi )[2u(\xi ) -v(\xi )] d\xi \label{C}\\  
 v(x) &=&  \sin (\omega _n^-x) +\frac{1}{\omega _n^-} \int _0^x 
  \sin (\omega _n^-(x-\xi )) \phi ^2(\xi )[u(\xi ) -2v(\xi )] d\xi \label{D}  
\end{eqnarray}
\end{subequations}
Standard calculations yield the explicit
formula for the Green's function $\Gamma _{\omega _n^+}$, which is given by
\begin{eqnarray}
\left\{ \begin{array}{l} 
\Gamma _{\omega _n^+}(x,\xi) = 
 \displaystyle \frac{1}{4\omega _n^+ \sinh(\omega _n^+ )} \Big\{ [\sinh(\omega _n^+\xi)e^{-\omega _n^+}-\sinh(\omega _n^+(1-\xi))]e^{\omega _n^+ x} \\
 \qquad \qquad \qquad \qquad \qquad \quad  -[\sinh(\omega _n^+\xi)e^{\omega _n^+}
 -\sinh(\omega _n^+ (1-\xi))]e^{-\omega _n^+ x}\Big\} \\
  \qquad \qquad \quad +\displaystyle \frac{\sinh(\omega _n^+ |x-\xi |)}{2\omega _n^+ } \end{array}
  \right. \label{Gamma1} 
\end{eqnarray}   
The function $|\Gamma _{\omega _n^+} (x,\xi )|$ 
assumes its maximum on $[0,1]^2$ at the point 
$(x,\xi ) = \left(\frac{1}{2} , \frac{1}{2} \right) $ and its  
maximum value is given by 
\begin{equation}
 \left|\Gamma _{\omega _n^+} \left({\textstyle\frac{1}{2},\frac{1}{2} }\right)\right| =
\frac{\sinh ^2(\frac{\omega _n^+ }{2})}{\omega _n^+ \sinh(\omega _n^+ )} = \frac{\cosh(\omega _n^+ )-1}{2\omega _n^+ \sinh(\omega _n^+ )} = \mathcal{O}\left( \frac{1}{\omega _n^+ }\right) \quad 
(\omega _n^+ \to \infty ) \label{omega}
\end{equation}
(uniformly in $\omega _n^+ $). Now it is a matter of routine
estimates (combining the representation
(\ref{C}),(\ref{D})  with the uniform 
estimate (\ref{omega}) and the property 
$\omega _n^\pm \sim n $, as $n\to \infty $)
to verify the assertion of the lemma. \hfill $\Box $    
\subsection{Verification of the properties listed in Section \ref{Spec}}
\label{Ver} 
\begin{enumerate}
 \item[\ref{1}.] Lemma \ref{LK} (i).
 \item[\ref{2}.] Let $\lambda \in \C\setminus \{ 0\} $ be an eigenvalue
for $L$. Then the argument on page 13 of \cite{CGNT06} (``Claim'')
shows that $\lambda ^2$ is real.  We can then use the argument 
on pages 49 and 50 of \cite{RSS05} to show that $\lambda $ itself is 
purely imaginary. This argument uses, besides the properties of 
$L_+$ listed in Lemma \ref{L2}, the convexity property
 $\frac{d}{d\muz } \langle    \phi ,\phi \rangle      = 2 \langle    \partial _{\muz } \phi ,\phi \rangle      \, >0$ (Lemma \ref{L5}). 
 \item[\ref{3}.] The first assertion (all but finitely many eigenvalues are 
simple) was already mentioned in Lemma \ref{LK} (i). \\[0.1ex]
A proof of the second one (there are no generalized eigenvectors
associated with non-zero eigenvalues) can be found in  
\cite[pages 50-51]{RSS05}. The proof refers to the whole-space 
case, but carries over to the zero-boundary case if 
the domains of the various operators 
involved are modified suitably. \\[0.1ex]
The third assertion (each eigenspace is at most two-dimensional)
may be seen as follows. Let $\lambda \in \C $ be an eigenvalue and 
$\textrm{Eig}(\lambda )$ the corresponding eigenspace and consider
the linear map $F:\textrm{Eig} (\lambda )\to \C ^2 $, defined by 
$FV := V'(0)$. Using the left boundary condition ($V(0)={0\choose 0}$) 
and ODE uniqueness, is is easy to see that $F$ is one-to-one.
Thus $\dim (\textrm{Eig}(\lambda ))=\dim (\textrm{im}(F)) \le 2$.   
  \item[\ref{4}.] Let ${u\choose v}$ be an eigenfunction for $L^*$ 
associated with the eigenvalue $\lambda =0$. Then
$L_-u=0$ and $L_+v=0$. From Lemma \ref{L1} (i) and Lemma \ref{L2}
 we get $u=c\phi $ ($c\in \C $) and $v=0$; hence 
\[ {u\choose v} \in \textrm{span} \left\{ {\phi \choose 0} \right\} .\] 
It is easy to see that $W_1 ={0 \choose \partial _{\muz } \phi  } $ 
satisfies  $-L^*W_1=V_1$; it is therefore a generalized eigenvector.
We want to show that there cannot be another linearly independent 
generalized eigenvector. To prove this, let's assume that 
 ${u\choose v}$ is such an ``additional'' generalized eigenvector, i.e. 
\[ -L^*{{u\choose v} = W_1 = {0 \choose \partial _{\muz } \phi  }} ,\] 
since the eigenspace is one-dimensional. In particular, 
\begin{equation}
L_-u=-\partial _{\muz } \phi .\label{gr}
\end{equation} 
By the Fredholm alternative, this implies
\begin{equation}
 0 = 2 \langle    \phi ,\partial _{\muz } \phi \rangle      = \frac{d}{d\muz } \langle    \phi ,\phi \rangle      
\label{fred}
\end{equation}
since $\phi $ is a non-trivial solution to $L_-\phi =0$.  But 
this contradicts Lemma \ref{L5}. As a result, no solution $u$ to 
(\ref{gr}) can exist, so there is no ``additional'' generalized 
eigenvector.      
\item[\ref{6}.] Basis property: Lemma \ref{LK} (i). \\[1ex] 
 The Bessel-sequence property follows from Lemma \ref{Lkam}: 
Since it is obviously sufficient to establish
 that the sequence of (generalized) eigenvectors
 for $M$ is a Bessel sequence, we will show that  
 \begin{equation}
  \forall W=\textstyle {u\choose v}\in [L^2(0,1)]^2: \qquad(a_n^\pm )_{n\in \N }
  :=\left( \langle    W,W_n^\pm \rangle     \right) _{n\ge \N }  \in \ell ^2.
  \label{bessel}
 \end{equation} 
 Clearly, we may skip a finite number of terms in (\ref{bessel}). 
 For simplicity, we will also restrict ourselves to $a_n^+$.  
Let $n_0$ be an index such that 
 \[ W_n^+(x) = \sin (\omega _nx) {0\choose 1} + \frac{R_n^+(x)}{n} \quad (\forall x\in [0,1],\, n\ge n_0),\]
 where $\omega _n:=\sqrt{\beta _n-\muz } $ and $R_n^+:[0,1]\to \C^2$ are 
 continuous functions satisfying $|R_n^+(x)|\le C$ uniformly in $x$ and $n$
 (see Lemma \ref{Lkam}).  Now
 \[ a_n^+ =\int _0^1 v(x) \sin (\omega _nx) dx + \frac{\langle    W,R_n^+\rangle      }{n} =: b_n^+ +c_n^+ .\]  
 The sequence $(c_n^+)$ is easily seen to be an $\ell ^2$-sequence:
 \[ |\langle    W,R_n^+\rangle     |\le C \int _0^1|W(x)|dx \le C\|W\| 
\quad \Rightarrow \quad   |c_n^+| \le \frac{C\|W\|}{n} \in \ell ^2.\]   
 To see that $(b_n^+)$ is square-summable as well, we note that, 
by Lemma \ref{LK}, the sequence $(\omega _n)_{n\ge n_0} $ has the 
asymptotics $\omega _n\sim \pi n$, as $n\to \infty $, and is 
  therefore separated. By \cite[Theorem 3.4]{Y80}, this implies 
that the exponential system 
  $\{ e^{i\omega _nx}\} _{n\ge n_0}$ forms a Bessel sequence in $L^2(0,1)$. 
Thus 
  \[ b_{j,n}^+ := \int _0^1 v_j(x) e^{i\omega _nx} dx \in \ell ^2  \qquad (j\in \{1,2\} ), \]
  where $v_1(x) := \Re (v(x))$ and $v_2(x):=\Im (v(x))$, and so
  \[ b_n^+ = \Im (b_{1,n}^+) +i \Im (b_{2,n}^+) \in \ell ^2 ,\]
  which completes the proof of $(a_n^+)\in \ell ^2$.        
\item[\ref{7}.] This also follows from Lemma \ref{Lkam}. 
Let $\Omega '=(a,b)\subset (0,1)$. Then we have
\begin{eqnarray} 
\int _{\Omega '}\sin^2(\omega x)dx
  \longrightarrow \frac{b-a}{2} ,
\quad \textrm{as $\omega \to \infty $. }\label{lim} \nonumber
 \end{eqnarray}
As above, we  appeal to Lemma \ref{Lkam} to write $W_n^+ $ in the form
\[ W_n^+(x) = \sin (\omega _nx) {0\choose 1} + \frac{R_n^+(x)}{n} \quad
(\forall x\in [0,1],\, n\ge n_0).\]
Now
 \begin{eqnarray*} \int _{\Omega '}|W_n^+(x)|^2 dx &\ge & 
  \int _{\Omega '}\sin ^2(\omega _nx) dx - \frac{2}{n} 
   \int _{\Omega '}{\textstyle \left|\sin (\omega _nx)\left\langle     {0\choose 1},R_n^\pm
   (x)\right\rangle      \right|} dx - \cdots \\
 & & \qquad - \frac{1}{n^2} \int _{\Omega '}|R_n^\pm (x)|^2dx \\
   &\ge &   \int _{\Omega '}\sin ^2(\omega _nx) dx - \frac{2C(b-a)}{n} -\frac{C^2(b-a)}{n^2} 
   \longrightarrow \frac{b-a}{2} ,
   \end{eqnarray*} 
as $n\to \infty $ (note that $\lim _{n\to \infty } \omega _n = \infty $). 
This concludes the proof.  
\item[\ref{9}.] clear
\item[\ref{10}.] According to Lemma \ref{LK}, we may write 
\[ \beta _n = \pi ^2 n^2 +r_n  \qquad (n\ge n_0),\]
where $|r_n|\le C$ for some constant $C>0$ (independent of $n$).  
It follows that  
\[ \beta _{n+1} -\beta _n = \pi ^2[(n+1)^2-n^2] + r_{n+1}-r_n 
= \pi ^2[2n+1] + r_{n+1}-r_n \to \infty ,\]
as $n\to \infty $, which, by \cite[Corollary]{Y98}, implies that 
the sequence 
\[ \mathscr{E} :=\{ e^{-i\beta _nt},1, e^{i\beta _nt} \}_{n\ge 2}\]
is a Riesz-Fischer sequence in $L^2(0,T)$ for every $T>0$. 
We will show that adding
the function $f(t)=t$ to $\mathscr{E} $ will result in a Riesz-Fischer 
sequence in $L^2(0,T)$ as well, \textit{if $T>0$ is large enough.}   
For ease of notation, let's define 
\[ e_0(t) = 1, \quad e_{\pm m}(t) := e^{\pm i\beta _{m+1}t} \quad (m\ge 1) .\] 
The fact that the sequence $\mathscr{E} $ is a Riesz-Fischer means that
there is a constant $A _{\mathscr{E}}>0$ such that 
\begin{equation}
\int _0^T \left| \sum _{n=-\infty }^\infty a_ne_n (t) \right| ^2 
dt \ge A _{\mathscr{E}}\sum _{n=-\infty }^\infty |a_n|^2 
\qquad (\forall (a_n)_{n\in \Z } \in \ell ^2)
\label{RF1}    
\end{equation} 
Thus
\begin{eqnarray*}
\int _0^T \left|\sum _{n=-\infty }^\infty a_ne_n (t) \, + bf(t)  \right| ^2 
dt 
 &\stackrel{(\ref{RF1})}{\ge } & A _{\mathscr{E}} \| a\| _{\ell ^2}^2
  + \| f\| _{L^2(0,T)}^2 |b|^2 -
 2 |b| \, \| a\| _{\ell ^2}  \| F\| _{\ell ^2}  
\end{eqnarray*} 
where the sequence $(F_n)_{n\in \Z }$ is defined by 
\[ F_{\pm m} := \int _0^T e_{\pm m}(t) f(t) dt = 
\int _0^T e^{\pm i\beta _{m+1}t}f(t) dt \quad (m\ge 0).\]   
(Since the sequence $(-\beta _m,\beta _m )_{m\in \N }$ is obviously  
separated, the sequence $\mathscr{E} $ is a Bessel sequence, which 
implies  that $F\in \ell ^2$. Since $f(t)=t$, the 
 $F_n$'s can also be computed explicitly to verify the square-summability 
of $F$.) Utilizing the elementary Young's inequality, we can   
continue the estimation above to obtain  
\begin{eqnarray*}
\int _0^T \left|\sum _{n=-\infty }^\infty a_ne_n (t) \, + bf(t)  \right| ^2 
dt &\ge & \varepsilon  \| a\| _{\ell ^2}^2
  + \left(\| f\| _{L^2(0,T)}^2 - \frac{\| F\| _{\ell
        ^2}^2}{A_{\mathscr{E}}-\varepsilon } \right) |b|^2   
\end{eqnarray*}  
(for $\varepsilon \in (0,A_{\mathscr{E}})$), which 
will yield the assertion provided that the  condition 
\begin{equation}
 \| F\| _{\ell ^2}^2 <(A_{\mathscr{E}} -\varepsilon ) \| f\| _{L^2(0,T)}^2 .
\label{RF2}
\end{equation}   
is satisfied. Clearly, $ \| f\| _{L^2(0,T)}^2 =\frac{T^3}{3} $,
 since $f(t) =t$. Moreover, it is a matter of routine calculations to verify
that $\| F\| _{\ell ^2}^2 \le C_fT^2$  
for some $T$-independent constant $C_f$. 
Thus, condition (\ref{RF2}) takes the form
 \begin{equation}
 C_f <\frac{T}{3} (A_{\mathscr{E}} -\varepsilon )  
\label{RF5}
\end{equation}
(note that the constant $ A_{\mathscr{E}}$ depends on $T$, 
$A_{\mathscr{E}}= A_{\mathscr{E}}(T)$). 
 Finally, it is easy to see that $A_{\mathscr{E}}$ can be chosen 
such that $A_{\mathscr{E}}(T_2)\ge A_{\mathscr{E}}(T_1)$ if
$T_2\ge T_1$, which implies that condition (\ref{RF5})
(and hence (\ref{RF2})) can be 
fulfilled by choosing $T$ large enough.  This concludes the proof that
the sequence $\mathscr{E} \cup \{ f\} $ is a  
the Riesz-Fischer sequence in $L^2(0,T)$ if $T>0$ is sufficiently large.
\end{enumerate}
\textbf{Acknowledgment.} This research  was supported by
 the Natural Sciences and Engineering Research Council
of Canada through its Discovery Grant programme.
The authors would like to thank R.~Illner for valuable discussions. 

\end{document}